\documentclass[a4 paper,12pt]{article}
\usepackage[T1]{fontenc}
\usepackage[latin1]{inputenc}
\usepackage{hyperref}
\usepackage[dvips]{graphics}
\usepackage[english]{babel}
\usepackage{amsmath}
\usepackage{amssymb}
\usepackage{amsopn}
\usepackage{amsbsy}
\usepackage{amstext}
\usepackage{latexsym}
\usepackage{amsfonts}
\usepackage{array}
\usepackage{epsfig}
\usepackage{mathrsfs}
\title{An identity involving the least common multiple of binomial coefficients and its application}
\author{\sc Bakir FARHI}
\date{}

\newtheorem{thm}{Theorem}
\newtheorem{prop}[thm]{Proposition}

\newtheorem{coll}[thm]{Corollary}

\newtheorem{thmn}{Theorem$\!\!\!$}

\let\epsilon=\varepsilon
\def\lcm{{\rm lcm}}

\def\EMdash{\leavevmode\hbox to 7.5mm{\vrule height .63ex depth -.59ex
    width 5.4mm\hfill}}

\begin{document}
\maketitle \vspace{-7cm}
\begin{flushleft}
To appear in {\it American Mathematical Monthly}
\end{flushleft}~\vspace{4cm}
\begin{center}
{\tt bakir.farhi@gmail.com}
\end{center}
\begin{abstract}
In this paper, we prove the identity $$\lcm\left\{\binom{k}{0} ,
\binom{k}{1} , \dots , \binom{k}{k}\right\} = \frac{\lcm(1 , 2 ,
\dots , k , k + 1)}{k + 1} ~~~~~~ (\forall k \in \mathbb{N}) .$$
As an application, we give an easily proof of the well-known
nontrivial lower bound $\lcm(1 , 2 , \dots , k) \geq 2^{k - 1}$
$(\forall k \geq 1)$.
\end{abstract}
{\bf MSC:} 11A05. \\
{\bf Keywords:} Least common multiple; Binomial coefficients;
Kummer's theorem.

\section{Introduction and Results}

Many results concerning the least common multiple of a sequence of
integers are known. The most famous is nothing else than an
equivalent of the prime number theorem; it states that $\log\lcm(1
, 2 , \dots , n) \sim n$ as $n$ tends to infinity (see, e.g.,
\cite{hw}). Effective bounds for $\lcm(1 , 2 , \dots , n)$ are
also given by several authors. Among others, Nair \cite{n}
discovered a nice new proof for the well-known estimate $\lcm(1 ,
2 , \dots , n) \geq 2^{n - 1}$ $(\forall n \geq 1)$. Actually,
Nair's method simply exploits the integral $\int_{0}^{1} x^n (1 -
x)^n d x$. Further, Hanson \cite{ha} already obtained the upper
bound $\lcm(1 , 2 , \dots , n) \leq 3^n$ $(\forall n \geq 1)$.

Recently, many related questions and many generalizations of the
above results have been studied by several authors. The interested
reader is referred to \cite{b}, \cite{f1}, and \cite{hf}.

In this note, using Kummer's theorem on the $p$-adic valuation of
binomial coefficients (see, e.g., \cite{kummer}), we obtain an
explicit formula for $\lcm\{\binom{k}{0} , \binom{k}{1} , \dots ,
\binom{k}{k}\}$ in terms of the least common multiple of the first
$k + 1$ consecutive positive integers. Then, we show how the
well-known nontrivial lower bound $\lcm(1 , 2 , \dots , n) \geq
2^{n - 1}$ $(\forall n \geq 1)$ can be deduced very easily from
that formula. Our main result is the following:
\begin{thm}\label{t1}
For any $k \in \mathbb{N}$, we have:
$$\lcm\left\{\binom{k}{0} , \binom{k}{1} , \dots , \binom{k}{k}\right\} = \frac{\lcm(1 , 2 , \dots , k , k + 1)}{k + 1} .$$
\end{thm}

First, let us recall the so-called Kummer's theorem:
\begin{thmn}[Kummer \cite{kummer}]
Let $n$ and $k$ be natural numbers such that $n \geq k$ and let
$p$ be a prime number. Then the largest power of $p$ dividing
$\binom{n}{k}$ is given by the number of borrows required when
subtracting $k$ from $n$ in the base $p$.
\end{thmn}

Note that the last part of the theorem is also equivalently stated
as the number of carries when adding $k$ and $n - k$ in the base
$p$.

As usually, if $p$ is a prime number and $\ell \geq 1$ is an
integer, we let $v_p(\ell)$ denote the normalized $p$-adic
valuation of $\ell$; that is, the exponent of the largest power of
$p$ dividing $\ell$. We first prove the following proposition.

\begin{prop}\label{prop1}
Let $k$ be a natural number and $p$ a prime number. Let $k =
\sum_{i = 0}^{N} c_i p^i$ be the $p$-base expansion of $k$, where
$N \in \mathbb{N}$, $c_i \in \{0 , 1 , \dots , p - 1\}$ (for $i =
0 , 1 , \dots , N$) and $c_N \neq 0$. Then we have:
$$\max_{0 \leq \ell \leq k} v_p\left(\binom{k}{\ell}\right) = v_p\left(\binom{k}{p^N - 1}\right) = \begin{cases}0 & \text{if $k = p^{N + 1} - 1$}
\\ N - \min\{i ~|~ c_i \neq p - 1\} & \text{otherwise.}\end{cases}$$
\end{prop}

\noindent{\bf Proof.} We distinguish the following two cases:\\
{\bf 1\textsuperscript{st} case.} If $k = p^{N + 1} - 1$:\\
In this case, we have $c_i = p - 1$ for all $i \in \{0 , 1 , \dots
, N\}$. So it is clear that in base $p$, the subtraction of any
$\ell \in \{0 , 1 , \dots , k\}$ from $k$ doesn't require any
borrows. It follows from Kummer's theorem that
$v_p\left(\binom{k}{\ell}\right) = 0$, $\forall \ell \in \{0 , 1 ,
\dots , k\}$. Hence
$$\max_{0 \leq \ell \leq k}
v_p\left(\binom{k}{\ell}\right) = v_p\left(\binom{k}{p^N -
1}\right) = 0 ,$$ as required.\newpage

\noindent{\bf 2\textsuperscript{nd} case.} If $k \neq p^{N + 1} - 1$:\\
In this case, at least one of the digits of $k$, in base $p$, is
different from $p - 1$. So we can define:
$$i_0 := \min\{i ~|~ c_i \neq p - 1\} .$$
We have to show that for any $\ell \in \{0 , 1 , \dots , k\}$, we
have $v_p(\binom{k}{\ell}) \leq N - i_0$, and that
$v_p(\binom{k}{p^N - 1}) = N - i_0$. \\
Let $\ell \in \{0 , 1 , \dots , k\}$ be arbitrary. Since (by the
definition of $i_0$) $c_0 = c_1 = \dots = c_{i_0 - 1} = p - 1$,
during the process of subtraction of $\ell$ from $k$ in base $p$,
the first $i_0$ subtractions digit-by-digit don't require any
borrows. So the number of borrows required in the subtraction of
$\ell$ from $k$ in base $p$ is at most equal to $N - i_0$.
According to Kummer's theorem, this implies that
$v_p(\binom{k}{\ell}) \leq N - i_0$.

\noindent Now, consider the special case $\ell = p^N - 1 = \sum_{i
= 0}^{N - 1} (p - 1) p^i$. Since $c_0 = c_1 = \dots = c_{i_0 - 1}
= p - 1$ and $c_{i_0} < p - 1$, during the process of subtraction
of $\ell$ from $k$ in base $p$, each of the subtractions
digit-by-digit from the rank $i_0$ to the rank $N - 1$ requires a
borrow. It follows from Kummer's theorem that $v_p(\binom{k}{p^N -
1}) = N - i_0$. This completes the proof of the
proposition.\penalty-20\null\hfill$\blacksquare$\par\medbreak

Now we are ready to prove our main result.

\noindent{\bf Proof of Theorem \ref{t1}.} The identity of Theorem
\ref{t1} is satisfied for $k = 0$. For the following, suppose $k
\geq 1$. Equivalently, we have to show that
\begin{equation}\label{eq1}
v_p\left(\lcm\left\{\binom{k}{0} , \binom{k}{1} , \dots ,
\binom{k}{k}\right\}\right) = v_p\left(\frac{\lcm(1 , 2 , \dots ,
k , k + 1)}{k + 1}\right) ,
\end{equation}
for any prime number $p$.\\
Let $p$ be an arbitrary prime number and $k = \sum_{i = 0}^{N} c_i
p^i$ be the $p$-base expansion of $k$ (where $N \in \mathbb{N}$,
$c_i \in \{0 , 1 , \dots , p - 1\}$ for $i = 0 , 1 , \dots , N$,
and $c_N \neq 0$). By Proposition \ref{prop1}, we have
\begin{equation}\label{eq2}
\begin{split}
v_p\left(\lcm\left\{\binom{k}{0} , \binom{k}{1} , \dots ,
\binom{k}{k}\right\}\right) &= \max_{0 \leq \ell \leq k}
v_p\left(\binom{k}{\ell}\right) \\
&\!\!\!\!\!\!\!\!\!\!\!\!\!\!\!\!\!\!\!\!\!\!\!\!\!\!\!\!\!\!\!\!=
\begin{cases}0 & \text{if $k = p^{N + 1} - 1$} \\ N - \min\{i ~|~
c_i \neq p - 1\} & \text{otherwise.}\end{cases}
\end{split}
\end{equation}
Next, it is clear that $v_p(\lcm(1 , 2 , \dots , k , k + 1))$ is
equal to the exponent of the largest power of $p$ not exceeding $k
+ 1$. Since (according to the expansion of $k$ in base $p$) the
largest power of $p$ not exceeding $k$ is $p^N$, the largest power
of $p$ not exceeding $k + 1$ is equal to $p^{N + 1}$ if $k + 1 =
p^{N + 1}$ and equal to $p^N$ if $k + 1 \neq p^{N + 1}$. Hence, we
have
\begin{equation}\label{eq3}
v_p\left(\lcm(1 , 2 , \dots , k , k + 1)\right) =
\begin{cases}N + 1 & \text{if $k = p^{N + 1} - 1$} \\ N &
\text{otherwise.}\end{cases}
\end{equation}
Further, it is easy to verify that
\begin{equation}\label{eq4}
v_p(k + 1) = \begin{cases}N + 1 & \text{if $k = p^{N + 1} - 1$}
\\ \min\{i ~|~ c_i \neq p - 1\} & \text{otherwise.}\end{cases}
\end{equation}
By subtracting the relation (\ref{eq4}) from the relation
(\ref{eq3}) and using an elementary property of the $p$-adic
valuation, we obtain
\begin{equation}\label{eq5}
v_p\left(\frac{\lcm(1 , 2 , \dots , k , k + 1)}{k + 1}\right) =
\begin{cases}
0 & \text{if $k = p^{N + 1} - 1$} \\ N - \min\{i ~|~ c_i \neq p -
1\} & \text{otherwise.}
\end{cases}
\end{equation}
The required equality (\ref{eq1}) follows by comparing the two
relations (\ref{eq2}) and
(\ref{eq5}).\penalty-20\null\hfill$\blacksquare$\par\medbreak

\section{Application to prove a nontrivial lower bound for $\lcm(1 , 2 , \dots , n)$}

We now apply Theorem \ref{t1} to obtain a nontrivial lower bound
for the numbers $\lcm(1 , 2 , \dots , n)$ $(n \geq 1)$.
\begin{coll}\label{c1}
For all integer $n \geq 1$, we have:
$$\lcm(1 , 2 , \dots , n) \geq 2^{n - 1} .$$
\end{coll}

\noindent{\bf Proof.} Let $n \geq 1$ be an integer. By applying
Theorem \ref{t1} for $k = n - 1$, we have:
\begin{eqnarray*}
\lcm(1 , 2 , \dots , n) & = & n \cdot \lcm\left\{\binom{n - 1}{0}
,
\binom{n - 1}{1} , \dots , \binom{n - 1}{n - 1}\right\} \\
& \geq & n \cdot \max_{0 \leq i \leq n - 1} \binom{n - 1}{i} \\
& \geq & \sum_{i = 0}^{n - 1}\binom{n - 1}{i} = 2^{n - 1} ,
\end{eqnarray*}
as required. The corollary is
proved.\penalty-20\null\hfill$\blacksquare$\par\medbreak

~\\


\begin{thebibliography}{9}
\bibliographystyle{plain}
\bibitem{b}
P. Bateman, J. Kalb, and A. Stenger, A limit involving least
common multiples, {\it Amer. Math. Monthly.} {\bf 109} (2002)
393-394.
\bibitem{f1}
B. Farhi, Nontrivial lower bounds for the least common multiple of
some finite sequences of integers, {\it J. Number Theory.} {\bf
125} (2007) 393-411.
\bibitem{ha}
D. Hanson, On the product of the primes, {\it Canad. Math. Bull.}
{\bf 15} (1972) 33-37.
\bibitem{hw}
G. H. Hardy and E. M. Wright, {\it The Theory of Numbers}, 5th
ed., Oxford University. Press, London, 1979.
\bibitem{hf}
S. Hong and W. Feng, Lower bounds for the least common multiple of
finite arithmetic progressions, {\it C. R. Math. Acad. Sci.
Paris.} {\bf 343} (2006) 695-698.
\bibitem{kummer}
E. E. Kummer, \"Uber die Erg\"anzungss\"atze zu den allgemeinen
Reciprocit\"atsgesetzen, {\it J. Reine Angew. Math.} {\bf 44}
(1852) 93-146.
\bibitem{n}
M. Nair, On Chebyshev-type inequalities for primes, {\it Amer.
Math. Monthly.} {\bf 89} (1982) 126-129.
\end{thebibliography}
\end{document}